\documentclass{article}

\usepackage{amsmath,amsfonts,amssymb}
\usepackage{mathdots}
\usepackage{graphicx}
\usepackage[colorlinks,linkcolor=blue,citecolor=blue]{hyperref}
\usepackage{color}

\newtheorem{theorem}{Theorem}[section]
\newtheorem{proposition}{Proposition}[section]

\begin{document}

\title{An interlacing result for Hermitian matrices in Minkowski space}

\author{
D.B. Janse van Rensburg\footnotemark[1],
A.C.M. Ran\footnotemark[2],
M. van Straaten\footnotemark[1]
}
\renewcommand{\thefootnote}{\fnsymbol{footnote}}
\footnotetext[1]{School~of~Mathematical and Statistical~Sciences,
North-West~University,
Research Focus: Pure and Applied Analytics,
Private~Bag~X6001,
Potchefstroom~2520,
South Africa.
E-mail: \texttt{
dawie.jansevanrensburg@nwu.ac.za, madelein.vanstraaten@nwu.ac.za}
}
\footnotetext[2]{Department of Mathematics,Vrije Universiteit Amsterdam, De Boelelaan
    1111, 1081 HV Amsterdam, The Netherlands
    and Research Focus: Pure and Applied Analytics, North-West~University,
Potchefstroom,
South Africa. E-mail:
    \texttt{a.c.m.ran@vu.nl}}
\date{}




\maketitle

\begin{abstract} In this paper we will look at the well known interlacing problem, but here we consider the result for Hermitian matrices in the Minkowski space, an indefinite inner product space with one negative square. More specific, we consider the $n\times n$ matrix $A=\begin{bmatrix} J & u\\ -u^* & a\end{bmatrix}$ with $a\in\mathbb{R}$, $J=J^*$ and $u\in\mathbb{C}^{n-1}$. Then $A$ is $H$-selfadjoint with respect to the matrix $H=I_{n-1}\oplus(-1)$. The canonical form for the pair $(A,H)$ plays an important role and the sign characteristic coupled to the pair is also discussed.
\end{abstract}




\section{Introduction}

The Hermitian matrix $H=I_{n-1} \oplus (-1)$ defines an indefinite inner product space with one negative square, this is sometimes called the Minkowski space. The formula  $[x,y]=\langle Hx,y\rangle$, with $x,y\in\mathbb{C}^n$ and where $\langle \cdot\, ,\cdot \rangle$ denotes the standard inner product, defines an indefinite inner product on $\mathbb{C}^n$. The function  $[x,y]$ satisfies all the properties of the standard inner product, with the exception that  $[x,x]$ may be nonpositive for $n\neq 0$. Basic elements of the theory of indefinite inner product spaces are summarised in \cite{GLR}. 

An $n\times n$ matrix $A$ is called $H$-selfadjoint if it is selfadjoint in the indefinite inner product given by $H$, or equivalently, if $HA=A^*H$. The spectrum $\sigma(A)$ of an $H$-selfadjoint matrix $A$ is symmetric relative to the real axis. The sizes of the Jordan blocks in the Jordan normal form of $A$ corresponding to eigenvalue $\lambda$ are equal to the sizes of the Jordan blocks corresponding to eigenvalue $\bar{\lambda}$, see Proposition 4.2.3 in \cite{GLR}. 

Canonical forms exist for pairs of matrices $(A,H)$ where $A$ is $H$-selfadjoint, see for example \cite{BMRRR}. However, for the specific pair of matrices discussed in this paper, namely 
$$
A=\begin{bmatrix} J & u\\ -u^* & a\end{bmatrix} \quad \text{and} \quad H=I_{n-1}\oplus(-1),
$$
the possible canonical forms are restricted to the ones below, see Section 5.6 in \cite{GLR}. There exists an invertible matrix $S$ such that $S^{-1}AS=A_1\oplus \cdots \oplus A_k$ and $S^*HS=H_1\oplus \cdots \oplus H_k$,
where the blocks in the canonical form are one of the following types:
\begin{enumerate}
\item $A_j=\lambda$, $H_j=\pm 1$, with $\lambda\in \mathbb{R}$, and only for one real eigenvalue the sign can be negative,
\item $A_j=\begin{bmatrix} a+bi & 0 \\ 0 & a-bi\end{bmatrix}$, $H_j=\begin{bmatrix} 0 & 1 \\1 & 0 \end{bmatrix}$, with $a\in\mathbb{R}$ and $b>0$, 
\item $A_j=J_2(\lambda)$, $H_j=\pm \begin{bmatrix} 0 & 1 \\1 & 0 \end{bmatrix}$, with $\lambda\in\mathbb{R}$,
\item $A_j=J_3(\lambda)$, $H_j=\begin{bmatrix} 0 & 0 & 1 \\ 0 & 1 & 0 \\ 1 & 0 & 0 \end{bmatrix}$, with $\lambda\in\mathbb{R}$.
\end{enumerate}
In addition, only one block of the forms in either 2, 3, or 4 can occur.

Note that for the specific matrices $A$ and $H$, we have $J^*=J$ and $a\in\mathbb{R}$. The goal of this paper is to find the relationship between the eigenvalues of the matrix $A$ and the matrix $J$ and how they are interlaced.

Interlacing problems for Hermitian matrices in a definite inner product space are well known, and results can be found in Section 4.3 in \cite{HJ}. Theorem 4.3.17 in \cite{HJ} gives Cauchy's interlacing theorem for a bordered Hermitian matrix. For an application to graphs and subgraphs, see for example the paper by Haemers \cite{H} and references there. 

Returning to the indefinite case, a different point of view is the inverse eigenvalue problem, and more precisely the periodic Jacobi inverse eigenvalue problem. See for example the paper by Xu, Bebiano and Chen, \cite{Bebianoetal} and references mentioned there. The paper \cite{Bebianoetal} is concerned with the reconstruction of a Jacobi matrix. We were inspired by some results in this paper and we were curious about the sign associated to the specific eigenvalues. We therefore explore this, applied to general $H$-selfadjoint matrices with one negative square.

In our paper we are concerned with how the eigenvalues of the selfadjoint matrix $J$ interlace with those of the $H$-selfadjoint matrix $A$, and with the sign corresponding to the eigenvalues of $A$ in the canonical form of the pair $(A,H)$. The proof for the interlacing of the eigenvalues follows the same line of argument as Lemma 3.3 in \cite{Bebianoetal}. 

We consider the characteristic polynomial of $A$. If $\lambda \notin \sigma(J)$, we use a standard argument to write $\lambda I-A$ as
$$
\lambda I-A=\begin{bmatrix}
\lambda I-J & -u\\ u^* & \lambda-a
\end{bmatrix}=\begin{bmatrix}
\lambda I-J&0\\
u^*& 1
\end{bmatrix}\begin{bmatrix}
I&-(\lambda I-J)^{-1}u\\
0&\lambda-a+u^*(\lambda I-J)^{-1}u
\end{bmatrix}.
$$
Hence the characteristic polynomial of $A$ becomes
\begin{equation}\label{eq:detlambda-A}
p(\lambda) = \det(\lambda I-A)=\det(\lambda I-J)\cdot (\lambda-a+u^*(\lambda I-J)^{-1}u).
\end{equation}
A part of  the second term in \eqref{eq:detlambda-A} is a realization of the general form
$$
D + C(\lambda I_n - A)^{-1}B.
$$
This realization will be minimal if $n$ is as small as possible, or equivalently, the pair $(A,C)$ is observable and the pair $(A,B)$ is controllable. For our case this reduces to the observability of the pair $(J,u^*)$.

Recall, see Theorem 3.2.1 in \cite{HRS}, that a pair $(A,C)$ with $C$ an $m\times n$ matrix and $A$ an $n\times n$ matrix is 
observable if 
$$\mathcal{N}:=\textup{Ker} \begin{bmatrix}
C\\CA\\CA^2\\\vdots\\CA^n
\end{bmatrix}=\{0\}.$$
The Hautus test for observability, see Theorem 3.2.2 in \cite{HRS}, states the following: the matrix pair $(A,C)$ is observable if and only if 
$$\textup{rank}\begin{bmatrix}
\lambda I-A\\C
\end{bmatrix}=n\,\,\, \textrm{for all }\,\, \lambda\in\sigma(A).$$
The Hautus test for the pair $(J,u^*)$ insures that all eigenvalues of $J$ appear as the poles of the expression $u^*(\lambda I - J)^{-1}u$, multiplicities included.\\

Section 2 of this paper describes the relationships between observability and the eigenvalues and Section 3 is concerned with the sign connected to the eigenvalues in the pair $(A,H)$. Finally, in the last section a few examples are given to verify and clarify some of the theory.

\section{Observability}

In this section we will prove several results following from the observability of the matrix pair $(J,u^*)$.

\begin{proposition}\label{prop:observe}
Let $J=J^*\in\mathbb{C}^{n-1 \times n-1}$ and $u\in\mathbb{C}^{n-1}$. If the pair $(J,u^*)$ is observable, i.e., $\textup{rank}\begin{bmatrix} \lambda I-J\\u^* \end{bmatrix}=n-1$ then:
\begin{itemize}
\item[(i)] $J$ has $n-1$ distinct eigenvalues;
\item[(ii)] $\sigma(J)\cap\sigma(A)=\emptyset$;
\item[(iii)] $A$ is nonderogatory.
\end{itemize}
\end{proposition}

{\bf Proof.}
(i) This follows immediately from the Hautus test. Indeed, since $J=J^*$ it only needs to be shown that for no eigenvalue $\lambda_0$ the corresponding eigenspace ${\rm Ker\, }(\lambda_0 I-J)$ has dimension two or higher. However, if this would be the case, then $\begin{bmatrix} \lambda_0 I -J \\ u^*\end{bmatrix}$ can have rank $n-2$ at most, thus violating the assumption of observability.
\\

(ii) From \eqref{eq:detlambda-A} and the fact that $(\lambda I-J)$ is invertible ($\lambda \notin\sigma(J))$, we have
$$
\frac{\det(\lambda I-A)}{\det(\lambda I-J)}= \lambda-a+u^*(\lambda I-J)^{-1}u.
$$
Thus, if $\sigma(J)\cap\sigma(A)\neq\emptyset$, then for some $\lambda_0\in\sigma(J)\cap\sigma(A)$, there will be a $(\lambda-\lambda_0)$ cancelling on the left, leaving us with a polynomial of degree $n-2$ in the denominator. Hence, there are at most $n-2$ poles but on the right-hand side we know there should be exactly $n-1$ poles, by the observability. This contradiction proves (ii).\\

(iii) Let $\lambda_0\in\sigma (A)$, i.e., $Ax=\lambda_0 x$, with $x=\begin{bmatrix}x_1\\x_2\end{bmatrix}\neq 0$. For $A=\begin{bmatrix} J & u\\-u^* & a \end{bmatrix}$ one obtains
$$
(\lambda_0 I-J)x_1-ux_2=0\quad\textup{and}\quad u^*x_1+(\lambda_0-a)x_2=0.
$$
If $x_2=0$, we have from the first equation that $(\lambda_0 I-J)x_1=0$, and since $\lambda_0\notin\sigma(J)$ by (ii), $(\lambda_0 I-J)$ is invertible and hence it implies $x_1=0$. This is a contradiction since $x$ is an eigenvector of $A$ belonging to $\lambda_0$, hence $x_2\neq 0$. Thus, one can solve for $x_1$ in terms of $x_2$ from the first equation which means the eigenvector corresponding to $\lambda_0$ is determined by $x_2$. Therefore, the $\dim\textup{Ker}(\lambda_0 I -A)=1$ and hence $A$ is nonderogatory.

\hfill $\Box$

If $(J,u^*)$ is not observable, the problem can be reduced to a situation where observability is satisfied.

\begin{proposition}
Assume $(J,u^*)$ is not observable, and let $\mathcal{N}=\cap_{j=0}^{n-2} {\rm Ker\,}u^*J^j$ be the unobservable subspace. 
With respect to $\mathbb{C}^{n-1}=\mathcal{N}\oplus \mathcal{N}^\perp$, write $J=J_1\oplus J_2$ as well as $u=\begin{bmatrix} 0 & u_2^*\end{bmatrix}$. Then 
$\sigma(A)\cap \sigma(J)=\sigma(J_1)$, and 
$$
\sigma(A)=\sigma(J_1) \cup\sigma\left(\begin{bmatrix} J_2 & u_2 \\ -u_2 &a \end{bmatrix} \right).
$$
\end{proposition}

{\bf Proof.}
Using the Kalman decomposition into the unobservable space and its orthogonal complement we can reduce to a situation where observability is satisfied, as the pair $(J_2,u_2^*)$ is observable.  
Writing
$$
A=\begin{bmatrix} J_1 & 0 & 0 \\ 0 & J_2 & u_2 \\ 0 & -u_2^* & a\end{bmatrix}
$$
the statements in the proposition easily follow.

\hfill $\Box$

\section{Interlacing}

In this section the main result of this article is presente. It contains the way the eigenvalues of matrices $A$ and $J$ are interlacing with each other, together with the sign corresponding to the eigenvalues of the matrix $A$. From Proposition~ \ref{prop:observe} we have that the eigenvalues of $A$ are precisely the $n$ complex zeros of the function
$\lambda -a +u^*(\lambda I-J)^{-1}u$ and this function has $n-1$ real distinct poles. Denote the eigenvalues of $J$ by $\mu_{n-1} < \mu_{n-2} < \cdots <\mu_1$ and introduce 
$$
g(\lambda)=-u^*(\lambda I-J)^{-1} u.
$$
Since $J=J^*$, there is a unitary matrix $V$ such that $V^*JV=D={\rm diag\, }(\mu_j)_{j=1}^{n-1}$.
Hence, one can write 
$$
g(\lambda)=-u^*V^*(\lambda I-D)^{-1} Vu=-\sum_{j=1}^{n-1} \frac{d_j}{\lambda -\mu_j},
$$
where $d_j=((Vu)_j)^2 >0$. Note that $g^\prime (\lambda)=\sum_{j=1}^{n-1}\frac{d_j}{(\lambda-\mu_j)^2} $ is positive where it is defined, and that $\lim_{\lambda\to\pm\infty}g(\lambda)=0$. 
Finally, the eigenvalues of $A$, which we shall denote by $\lambda_j$, $j=1, \ldots , n$, are the solutions of $h(\lambda)-g(\lambda)=0$, where $h(\lambda)=\lambda -a$. 

\begin{theorem}
Let $A$ be an $H$-selfadjoint matrix with $H=I_{n-1} \oplus (-1)$ (Hermitian and invertible) and let $A=\begin{bmatrix}J & u \\ -u^* & a \end{bmatrix}$ with $a\in\mathbb{R}$, $u\in\mathbb{C}^{n-1}$ and $J=J^*$. If the pair $(J,u^*)$ is observable, then the conditions of Proposition \ref{prop:observe} hold. Furthermore, let $\mu_{n-1} < \mu_{n-2} < \cdots <\mu_1$ denote the $n-1$ distinct real eigenvalues of $J$ and let $\lambda_1,\hdots, \lambda_n$ be the eigenvalues of $A$. Then the eigenvalues of $A$ and $J$ interlace in the following possible ways coupled with the appropriate sign for $\varepsilon$.
\begin{itemize}
\item[{\bf 1a}] $\lambda_n  <\lambda_{n-1} <\mu_{n-1} <\lambda_{n-2} <\cdots < \lambda_1 <\mu_1$,\, where the sign $\varepsilon =-1$ is associated with the Jordan block of size $1$ for the eigenvalue $\lambda_n$;
\item[{\bf 1b}] $\lambda_n =\lambda_{n-1}  <\mu_{n-1} <\lambda_{n-2} <\cdots < \lambda_1 <\mu_1$,\, where the sign $\varepsilon=-1$ is associated with a Jordan block of size $2$ with eigenvalue $\lambda_n=\lambda_{n-1}$;

\item[{\bf 2\phantom{a}}] $\mu_{n-1} <\lambda_{n-2} <\cdots < \lambda_1 <\mu_1$, \textrm{and} $\lambda_n=\overline{\lambda_{n-1}}\notin\mathbb{R}$,

\item[{\bf 3a}] $\mu_{n-1} <\lambda_n <\mu_{n-2} <\cdots <\lambda_3<\mu_1<\lambda_2 <\lambda_1 $,\, where the sign  $\varepsilon=-1$,  is associated with the Jordan block of size $1$ for eigenvalue $\lambda_1$;
\item[{\bf 3b}] $\mu_{n-1} <\lambda_n <\mu_{n-2} <\cdots <\lambda_3<\mu_1<\lambda_2 =\lambda_1$, where the sign $\varepsilon=1$ is associated with a Jordan block of size $2$ with eigenvalue $\lambda_1=\lambda_2$;

\item[{\bf 4a}] $\mu_{n-1}<\lambda_n <\mu_{n-2} <\cdots <\mu_{j+1} <\lambda_{j+2}<\lambda_{j+1}<\lambda_j<\mu_j<\cdots <\mu_2<\lambda_1<\mu_1$, where the sign $\varepsilon=-1$ is  associated with the Jordan block of size $1$ for eigenvalue $\lambda_{j+1}$;
\item[{\bf 4b}] $\mu_{n-1}<\lambda_n <\mu_{n-2} <\cdots <\mu_{j+1} <\lambda_{j+2}=\lambda_{j+1}<\lambda_j<\mu_j<\cdots <\mu_2<\lambda_1<\mu_1$, where the sign $\varepsilon=1$ is associated with the Jordan block of size $2$ with eigenvallue $\lambda_{j+2}=\lambda_{j+1}$;
\item[{\bf 4c}] $\mu_{n-1}<\lambda_n <\mu_{n-2} <\cdots <\mu_{j+1} <\lambda_{j+2}<\lambda_{j+1}=\lambda_{j}<\mu_j<\cdots <\mu_2<\lambda_1<\mu_1$, where the sign $\varepsilon=-1$ is  associated with the Jordan block of size $2$ with eigenvalue $\lambda_{j+1} = \lambda{j}$;
\item[{\bf 4d}] $\mu_{n-1}<\lambda_n <\mu_{n-2} <\cdots <\mu_{j+1} <\lambda_{j+2}=\lambda_{j+1}=\lambda_j<\mu_j<\cdots <\mu_2<\lambda_1<\mu_1$ where the sign $\varepsilon = 1$ is associated with a Jordan block of size $3$ with eigenvalue $\lambda_{j+2}=\lambda_{j+1}=\lambda_j$.
\end{itemize}
\end{theorem}

{\bf Proof.}
The interlacing of the eigenvalues of $A$ and $J$, i.e., the way $h(\lambda)$ and $g(\lambda)$ intersect, follows a similar result by Lemma 3.3 in \cite{Bebianoetal}. Because of the fact that $A$ is nonderogatory, in cases 1b, 3b, 4b and 4c there is a Jordan block of size two corresponding to the two eigenvalues that coincide, while in case 4d there is a Jordan block of size three corresponding to the three eigenvalues that coincide. See Section 5.6 in \cite{GLR}. \\
\indent We would like to know in the cases 1a, 3a and 4a in the list, which one of the eigenvalues has the negative sign in the sign characteristic, and in cases 1b, 3b, 4b and 4c what the sign corresponding to the Jordan block of size two is. In order to answer these questions we first recall one of the descriptions of the sign characteristic, see Section 5.1 in \cite{GLR}. First note that for every real $\lambda$, the matrix $\lambda H-HA$ is an $n\times n$ Hermitian  matrix, and hence has $n$ real eigenvalues, which are denoted by 
$\nu_1(\lambda), \ldots , \nu_n(\lambda)$. It turns out that these can be chosen to be analytic functions of the real variable $\lambda$, and this will be done for now. Let $\lambda_1, \ldots , \lambda_r$ be the real eigenvalues of $A$, and write for every $i=1, \ldots , r$ and $j=1, \ldots, n$ the function $\nu_j(\lambda)$ as 
$$
\nu_j(\lambda)=(\lambda-\lambda_i)^{m_{ij}}\rho_{ij}(\lambda),
$$
where $\rho_{ij}(\lambda_i)\not=0$ and is a real number. Then the nonzero numbers among $m_{i1}, \ldots , m_{in}$ are the sizes of the Jordan blocks of $A$ corresponding to $\lambda_i$, and the sign in the sign characteristic of the pair $(A,H)$ corresponding to the block of size $m_{ij}\not=0$ is the sign of the real number $\rho_{ij}(\lambda_i)$. In particular, 
\begin{center}
if $m_{ij}=1$, then $\rho_{ij}(\lambda_i)= \nu_j^\prime (\lambda_i)$, \\
\ \ if $m_{ij}=2$, then $\rho_{ij}(\lambda_i)=\frac{1}{2}\nu_j^{\prime\prime}(\lambda_i)$.
\end{center}

In order to find the signs in the particular situation we have at hand, we argue as follows. The fact that $\nu(\lambda)$ is an eigenvalue of $\lambda H-HA$ implies that
\begin{align*}
0&=\det(\nu I-(\lambda H-HA)) =\det \begin{bmatrix} \nu I-\lambda I +J & u \\ u^* & \nu +\lambda -a\end{bmatrix}\\
&=\det\left(
\begin{bmatrix} (\nu-\lambda)I+J & 0 \\ u^* & 1\end{bmatrix}\begin{bmatrix} I & ((\nu-\lambda)I+J)^{-1} u \\ 0 & (\nu+\lambda)-a -u^*((\nu-\lambda) I+J)^{-1} u\end{bmatrix} \right)\\
&=\det ((\nu-\lambda)I+J)\cdot \left( (\nu+\lambda)-a -u^*((\nu-\lambda) I+J)^{-1} u\right)\\
&=\det ((\nu-\lambda)I+J)\cdot \left( (\nu+\lambda)-a +u^*((\lambda-\nu) I-J)^{-1} u\right).
\end{align*}
It follows that $\nu$ satisfies $(\nu+\lambda)-a +u^*((\lambda-\nu) I-J)^{-1} u=0$, in other words,
$$
h(\nu+\lambda)-g(\lambda-\nu)=0,
$$
or more explicitly,
$$
\nu+\lambda-a=\sum_{j=1}^{n-1}\frac{d_j}{\lambda-\nu-\mu_j}.
$$
This determines $\nu$ implicitly as a function of $\lambda$. For fixed $\lambda$ we know that there have to be $n$ real solutions. Introduce 
$$
H(\lambda, \nu)=h(\lambda+\nu)-g(\lambda-\nu).
$$
When $m_{ij}\not=0$, we have $\nu_j(\lambda_i)=0$ and  $H(\lambda_i,0)=0$. Applying the implicit function theorem we obtain
$$
\nu_j^\prime(\lambda_i)=-\left({\frac{\partial H}{\partial \lambda}}/{\frac{\partial H}{\partial \nu}}\right)\rfloor_{\lambda=\lambda_i \atop \nu=0}.
$$
Now 
\begin{align*}
&\frac{\partial H}{\partial \lambda}\rfloor_{\lambda=\lambda_i \atop \nu=0} =
1-g^\prime(\lambda-\nu)\rfloor_{\lambda=\lambda_i \atop \nu=0} =1-g^\prime(\lambda_i),
\\
&\frac{\partial H}{\partial \nu}\rfloor_{\lambda=\lambda_i \atop \nu=0} =
1+g^\prime(\lambda-\nu)\rfloor_{\lambda=\lambda_i \atop \nu=0} =1+g^\prime(\lambda_i),
\end{align*}
so 
$$
\nu_j^\prime(\lambda_i) =-\left(\frac{1-g^\prime(\lambda_i)}{1+g^\prime(\lambda_i)}\right).
$$
Recall that $g^\prime(\lambda)>0$ whenever $\lambda$ is not one of the $\mu_j$'s.
Thus, if $m_{ij}=1$ the sign of $\nu_j^\prime(\lambda_i)$, and therefore also the sign attached to $\lambda_i$, is equal to the sign of $g^\prime(\lambda_i)-1$. We conclude that for an eigenvalue of multiplicity one,  the sign in the sign characteristic is determined by how $h$ and $g$ intersect at $\lambda_i$ as follows:
\begin{center}
if $g^\prime(\lambda_i) >1$ then the sign is $+1$ at eigenvalue $\lambda_i$; \\
if $g^\prime(\lambda_i)<1$ then the sign is $-1$ at eigenvalue $\lambda_i$. 
\end{center}

It remains to consider the signs of the Jordan blocks of order two. For that we first recall that $\nu(\lambda)$ satisfies $H(\lambda, \nu(\lambda))=0$. Taking the first derivative we have
$$
\frac{\partial H}{\partial \lambda}+\frac{\partial H}{\partial \nu}\cdot \nu^\prime(\lambda)=0,
$$
and differentiating this again gives
$$
\frac{\partial^2 H}{\partial \lambda^2}+\frac{\partial^2 H}{\partial \lambda\partial \nu}\cdot  \nu^\prime(\lambda)+\frac{\partial H}{\partial \nu}\cdot \nu^{\prime\prime}(\lambda)=0.
$$
A Jordan block of size two corresponds to a situation where $\nu_j(\lambda_i)=0$ and $\nu^\prime_j(\lambda_i)=0$, because $h$ and $g$ touch at $\lambda_i$. Solving the above equation for $\nu_j^{\prime\prime}(\lambda_i)$, we have
$$
\nu_j^{\prime\prime}(\lambda_i)=-\left(\frac{\partial^2 H}{\partial \lambda^2}/\frac{\partial H}{\partial \nu}\right)\rfloor_{\lambda=\lambda_i\atop \nu=0}=
-\frac{g^{\prime\prime}(\lambda_i)}{1+g^\prime(\lambda_i)}.
$$
Hence the sign corresponding to a block of size two with eigenvalue $\lambda_i$ is the sign of $-g^{\prime\prime}(\lambda_i)$. A little analysis shows that if the graph of $g$ locally around $\lambda_i$ lies above the graph of $h$, then the sign is $-1$ at eigenvalue $\lambda_i$, while if the graph of $g$ locally around $\lambda_i$ lies below the graph of $h$, then the sign is $+1$ at eigenvalue $\lambda_i$.
\hfill $\Box$

\section{Examples}

\begin{figure}[h!]
\includegraphics[height=7cm]{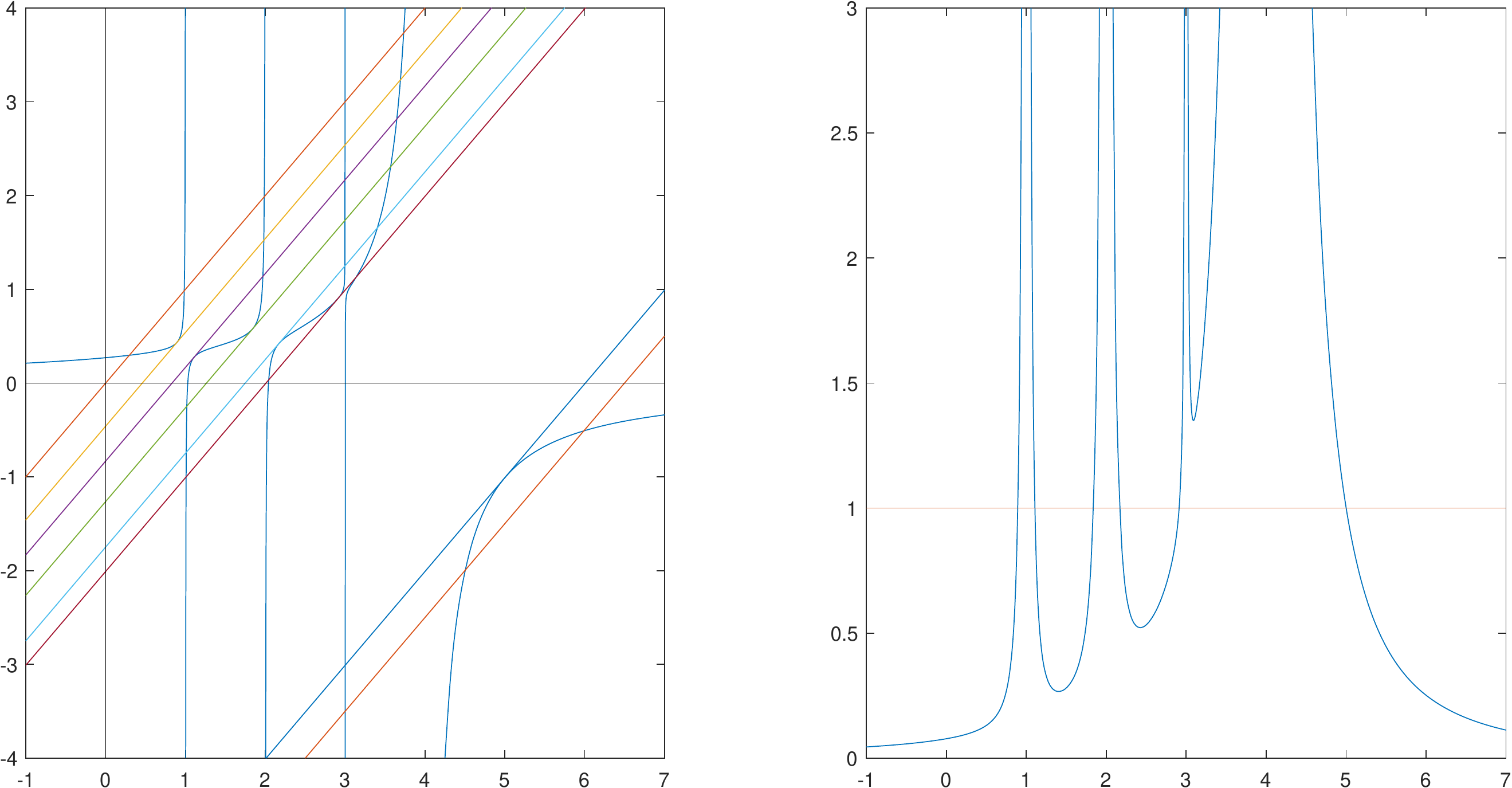}
\caption{On the left, some possible configurations of $h(\lambda)$ and $g(\lambda)$, on the right $g^\prime(\lambda)$ and the line $y=1$ to determine the points where there is a Jordan block of size two. }\label{Fig:1}
\end{figure}

\bigskip
As an example, consider the case where $J$ has eigenvalues $1, 2, 3, 4$, where the $d_j$'s are given by $0.01, 0.02, 0.001, 1$, respectively. Thus, 
$$
g(\lambda)= -\left(\frac{0.01}{\lambda-1}+\frac{0.02}{\lambda-2}+\frac{0.001}{\lambda-3}+\frac{1}{\lambda-4}\right).
$$
In Figure \ref{Fig:1}, the functions $g(\lambda)$ and $h(\lambda)=\lambda-a$ are plotted for several values of $a$, namely $a=0,\, 0.4591, \ 0.8319,\ 1.2631,\ 1.7485,\ 2.0087,\ 6.0097,\ 6.5$. 

We can draw the following conclusions regarding the sign of the eigenvalues depending on the choice of $a$:
\begin{enumerate}
\item In case 1a ($a=0)$ $\lambda_n$ comes with a negative sign, in case 3a ($a=6.5$) $\lambda_1$ comes with a negative sign, in case 4a (this would occur, e.g., for $a=1$) $\lambda_{j+1}$, the middle one of the three crossings in $(\mu_{j+1}, \mu_j)$ comes with a negative sign.
\item In case 1b ($a\approx 0.4591$) the sign corresponding to the block of order two with eigenvalue $\lambda_n=\lambda_{n-1}\approx 0.8934$ is $-1$, in case 3b ($a\approx 6.0097$) the sign corresponding to the block of order two with eigenvalue $\lambda_1=\lambda_2\approx 5,00155$ is $+1$, in case 4b ($a\approx 0.8319$ and $a\approx 1.7485$) the sign corresponding to the block of order two with eigenvalue $\lambda_{j+2}=\lambda_{j+1}\approx 1.10815$, respectively, $\lambda_{j+2}=\lambda_{j+1}\approx 2.1699$, is $+1$, and finally, in case 4c ($a\approx 1.2631$ and $a\approx 2.0087$) the sign corresponding to the block of order two with eigenvalue $\lambda_{j+1}=\lambda_j\approx 1.83895$, respectively, $\lambda_{j+1}=\lambda_j\approx 2.91185$, is $-1$. 
\end{enumerate}

In general, we can formulate the following when we view $a$ as a parameter. There is a value $a^-_1$ such that for $a<a^-_1$ we are in situation 1a, while for $a=a^-_1$ we are in situation 1b. Then there is a value $a^-_2$ such that for $a^-_1<a<a^-_2$ we are in case 2, while for $a=a^-_2$ we have either 4b for some $j$ between $n-2$ and $1$, or 4d, or 3b. From the other end, there is an $a^+_1$ such that for $a>a^+_1$ we are in the situation 3a, while for $a=a^+_1$ we are in the situation 3b. Finally, there is an $a^+_2$ such that for $a^+_2<a<a^+_1$ we are in the situation 2, while for $a=a^+_2$ we are either in the situation 4c,  or in the situation 4d, or in 1b. Eventually, for large positive $a$ there is one eigenvalue moving to $+\infty$, while there are 
$n-1$ eigenvalues approximating the eigenvalues of $J$ from the right. 

To check the latter statement, consider the equation $G(\lambda, a)=h(\lambda)-g(\lambda)=0$ as an equation determining $\lambda$ locally as a function of $a$. Then by the implicit function theorem
$$
\lambda_i^\prime(a)= -\left(\frac{\partial G}{\partial a}/\frac{\partial G}{\partial \lambda}\right) \rfloor_{\lambda=\lambda_i }= \frac{1}{ 1-g^\prime(\lambda_i)}.
$$
As for large values of $a$ the derivative $g^\prime(\lambda_i) >1$ (as can be seen from the graph of $g$) apart from the largest eigenvalue, we have $\lambda_i$ decreasing for $i=2,\ldots , n$ and increasing for $i=1$. Using the main result of \cite{RW} we obtain that the eigenvalues of $A$ as function of $a$ approximate the eigenvalues of $J$ as $a\to\infty$, with the exception of one, which goes to plus infinity.
\\

\begin{figure}[h!]
\phantom{.}\hskip -0.2cm
\includegraphics[height=5cm]{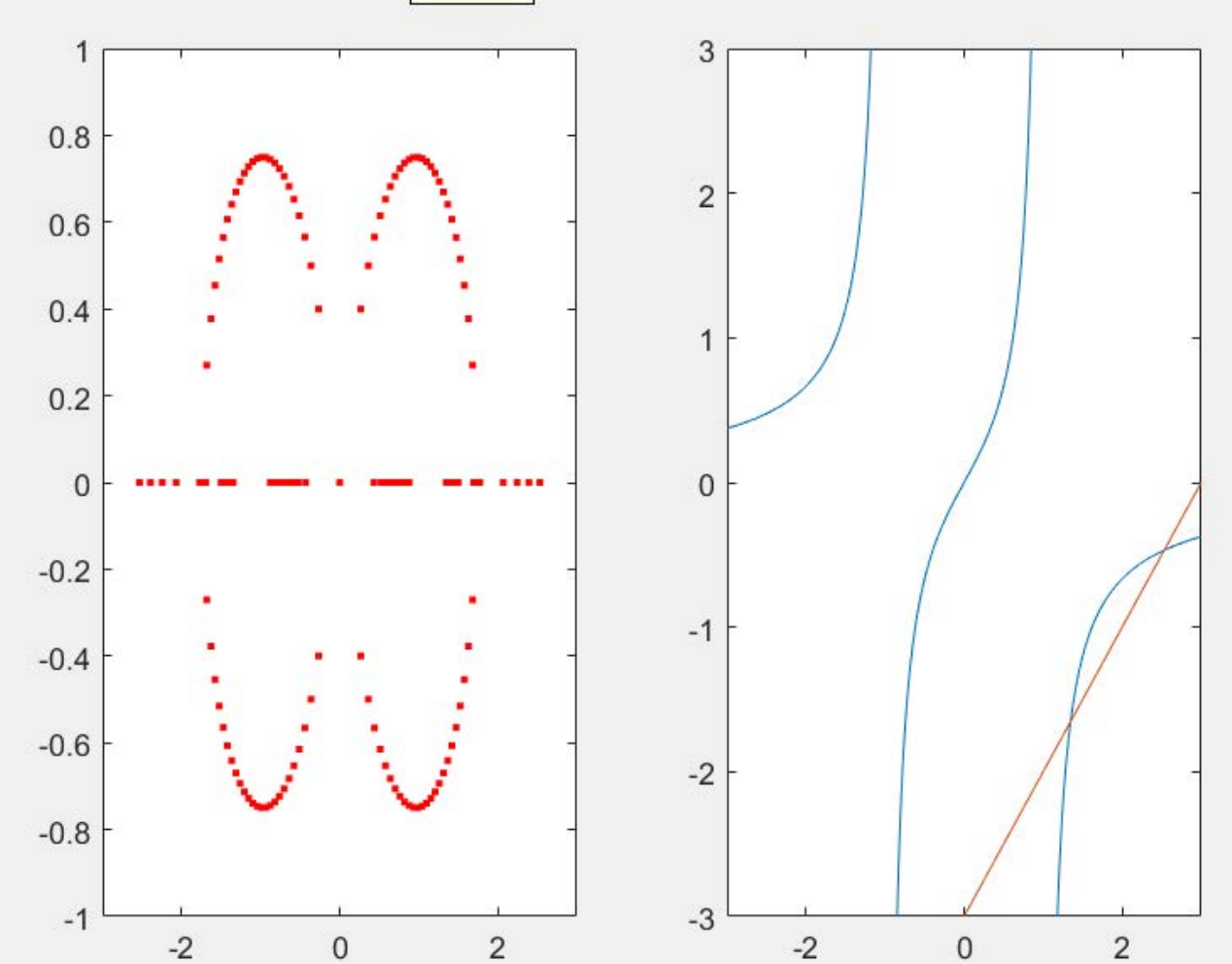}
\includegraphics[height=5cm]{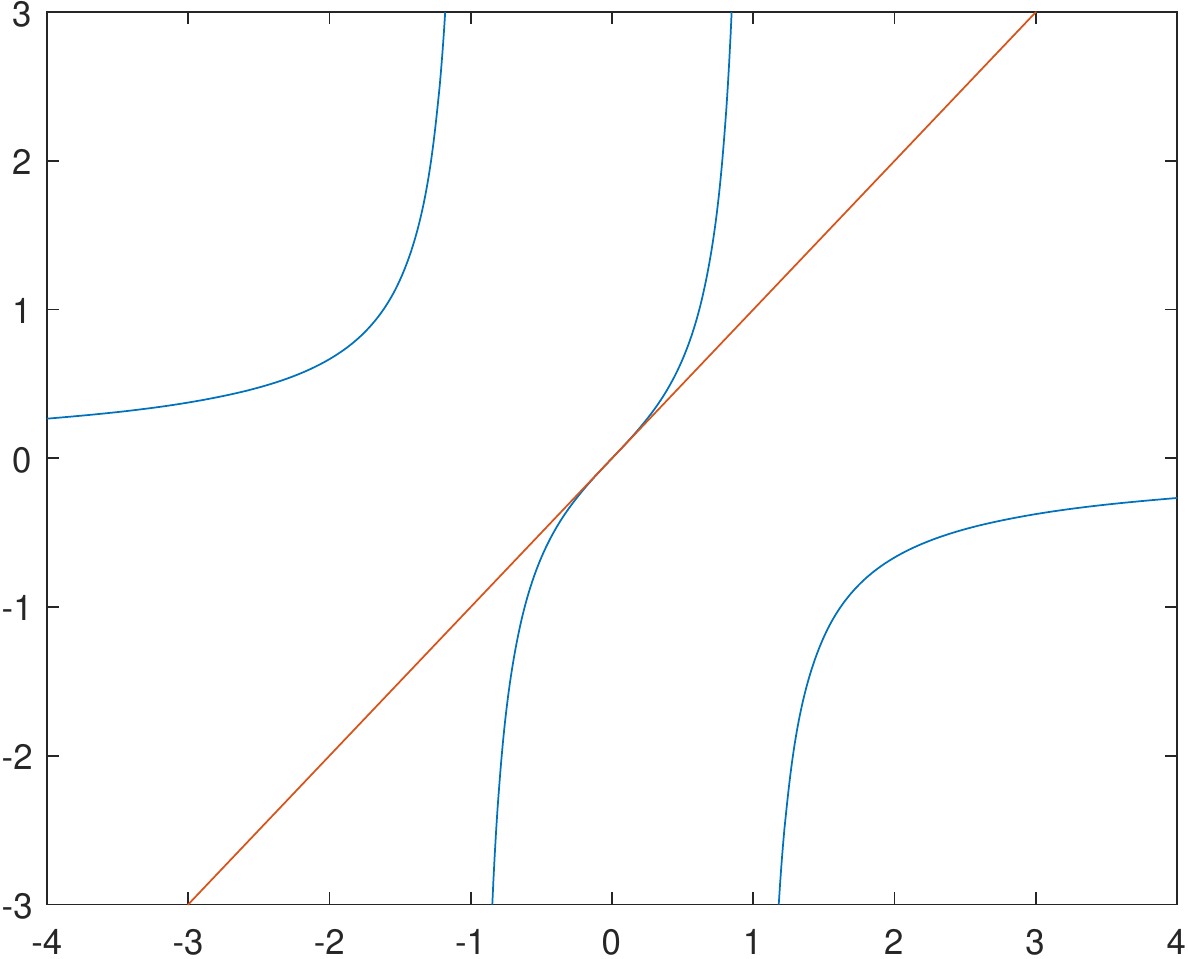}
\caption{The situation with a Jordan block of order three.}\label{Fig:2}
\end{figure}

The following example demonstrates a third order block when we take $\mu_1=1, \mu_2=-1$ and $u=[1/\sqrt{2}, 1/\sqrt{2}]$. In this case we have 
\[
A=\begin{bmatrix}1 & 0 &  1/\sqrt{2} \\ 0 & -1 & 1/\sqrt{2} 
\\ -1/\sqrt{2} & -1/\sqrt{2} & a \end{bmatrix}.
\] 
For $a=0$ there is a Jordan block of order three at the eigenvalue $0$.
In the left-hand graph in Figure \ref{Fig:2} the eigenvalues of $A$ are plotted for varying values of $a$. The dots in the circular-like shape are the complex eigenvalues that can occur, as the graph of $h(\lambda)$ moves from left to right over the curve of $g(\lambda)$ as function of $a$.

{\bf Acknowledgement.} 
The work of the first and third author is supported in part by the DSI-NRF Centre of Excellence in Mathematical and Statistical Sciences (CoE-MaSS, Grant Number 2022-012-ALG-ILAS).
The work of the second author is based on research supported in part by the National Research Foundation of South Africa (Grant Number 145688).

\end{document}